\documentclass[10pt]{amsart}
\usepackage{amsfonts}
\usepackage{amsmath}
\usepackage{amssymb}
\usepackage{amsmath,amssymb,amsfonts,amsthm,graphics,
latexsym, amscd, amsfonts, epsfig, eepic,epic}
\usepackage{mathrsfs}
\usepackage{algorithm}
\usepackage{algorithmic}
\usepackage{color}
\newcommand{\NEW}[1]{\begingroup\color{black}#1\endgroup}

\usepackage{eucal}
\usepackage{eufrak}
\usepackage[all]{xypic}
\usepackage{xspace}

\newtheorem{theorem}{Theorem}[section]

\newtheorem{lemma}[theorem]{Lemma}
\newtheorem{proposition}[theorem]{Proposition}
\newtheorem{corollary}[theorem]{Corollary}

\theoremstyle{remark}

\makeatletter \@addtoreset{equation}{section}

\title{Random 3-noncrossing partitions}

\author{Jing Qin  and Christian M. Reidys}
\email{reidys@santafe.edu}
\date{}
\begin{document}
\maketitle
\begin{abstract}
In this paper, we introduce polynomial time algorithms that generate
random $3$-noncrossing partitions and $2$-regular, $3$-noncrossing
partitions with uniform probability. A $3$-noncrossing partition
does not contain any three mutually crossing arcs in its canonical
representation and is $2$-regular if the latter does not contain
arcs of the form $(i,i+1)$. Using a bijection of Chen~{\it et al.}
\cite{Chen,Reidys:08tan}, we interpret $3$-noncrossing partitions
and $2$-regular, $3$-noncrossing partitions as restricted
generalized vacillating tableaux. Furthermore, we interpret the
tableaux as sampling paths of Markov-processes over shapes and
derive their transition probabilities.
\end{abstract}

\section{Introduction}\label{S:intro}
Recently, a paper written by Chen {\it et.~al.}~\cite{Chen:pnas}
attracts our attention. According to the bijection between the
$k$-noncrossing matchings and the oscillating tableaux \cite{Chen},
they identify the latter as stochastic processes over Young tableaux
of less than $k$ rows in order to uniformly generate a
$k$-noncrossing matching. Furthermore, since the generating function
of the corresponding oscillating lattice walks in $\mathbb{Z}^{k-1}$
that remains in the interior of the dominant Weyl chamber has been
given by Grabiner and Magyar \cite{Grabiner}, the key quantities,
the transition probabilities of the specific stochastic processes
can be derived with linear time complexity.

\begin{figure}[ht]
\centerline{%
\epsfig{file=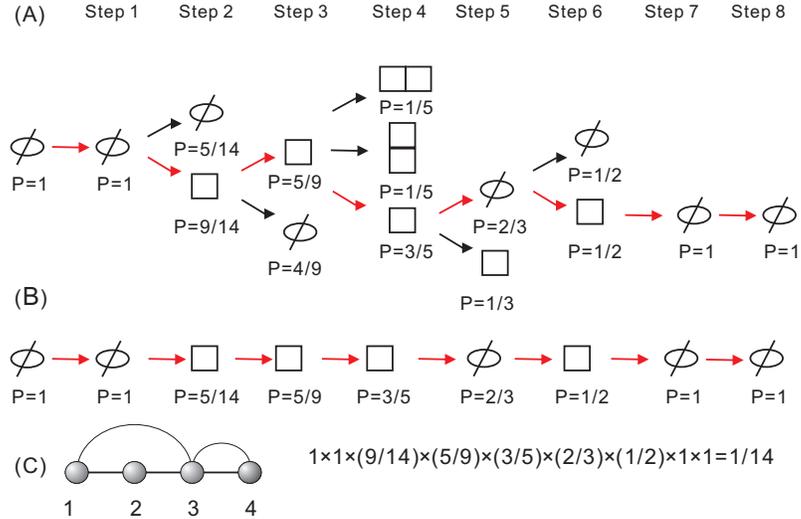,width=0.7\textwidth}\hskip15pt
 }
\caption{\small The idea behind the uniform generation: we consider
a stochastic process over shapes \textsf{(A)} extract a specific
sampling path \textsf{(B)} and \textsf{(C)} display the
corresponding $3$-noncrossing partition structure induced by this
path. The probabilities given in \textsf{(B)} are the conditional
probabilities with respect to their last step.  The transition
probabilities for given $k=3$ and $N$ are computed in Theorem
\ref{T:algpar} as a preprocessing step in polynomial time.}
\label{F:parcro}
\end{figure}

The objective to enumerate $k$-noncrossing partitions is much more
difficult since the corresponding lattice walks are not reflectable
\cite{Gessel:92} in $\mathbb{Z}^{k-1}$ in case of $k\geq 3$. Only
the case for $k=3$ has been solved by Bousquet-M\'{e}lou and Xin in
\cite{Xin:06} via the celebrated kernel method. Also in their paper,
they conjecture that $k$-noncrossing partitions are not
$P$-recursive for $k\geq 4$. For $k=3$, what Bousquet-M\'{e}lou and
Xin need in order to enumerate $3$-noncrossing partition is the
number of corresponding lattice walks starting and ending at
$(1,0)\in \mathbb{Z}^2$. However, our main idea is to interpret the
corresponding vacillating tableaux as sampling paths of
Markov-processes over shapes and derive their transition
probabilities, see Fig.~\ref{F:parcro}. In order to derive the
transition probabilities, what we need is the number of the
corresponding lattice walks ending at arbitrary $(i,j)\in
\mathbb{Z}^2$ such that $i>j\geq 0$.

This paper is organized as follows. Section~\ref{S:fact} describes
the basic facts of $3$-noncrossing partitions and $2$-regular,
$3$-noncrossing partitions. Section~\ref{S:vacillating} shows the
reader how we explore more information from the kernel equations. In
the meantime, we generate a $3$-noncrossing braid since we will
prove there exists a bijection between the set of $k$-noncrossing
partitions over $[N]= \{1,2,\dots, N\}$ and the set of
$k$-noncrossing braids over $[N-1]$. Furthermore, in the
Section~\ref{S:braid}, we will show the reader how to arrive at the
transition probabilities of corresponding lattice walks for
$2$-regular, $3$-noncrossing partitions from the fact that there
exists a bijection between the set of $2$-regular, $k$-noncrossing
partitions over $[N]$ and the set of $k$-noncrossing braids without
\NEW{loops} over $[N-1]$.

\section{Some basic facts}\label{S:fact}

\begin{figure}[ht]
\centerline{%
\epsfig{file=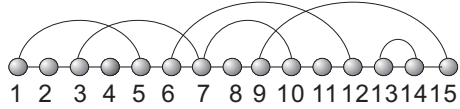,width=0.4\textwidth}\hskip15pt
 }
\caption{\small Standard representation and $k$-crossings: we
display the partition $\{2\}$, $\{8\}$, $\{11\}$, $\{1,5\}$,
$\{6,12\}$, $\{9,15\}$, $\{13,14\}$, $\{3,7,10\}$ of $[15]$.
Elements within blocks are connected in numerical order. The arcs
$\{(1,5),(3,7)\}$, $\{(3,7),(6,12)\}$ $\{(7,10),(9,15)\}$ and
$\{(6,12),(9,15)\}$ are all $2$-crossings. } \label{F:parcro1}
\end{figure}

A \emph{set partition} $P$ of $[N]$ is a collection of nonempty and
mutually disjoint subsets of $[N]$, called blocks, whose union is
$[N]$. A $k$-noncrossing partition is called \emph{$m$-regular},
$m\geq 1$, if for any two distinct elements $x$, $y$ in the same
block, we have $\vert x-y \vert \geq m$. A \emph{partial matching}
and a \emph{matching} is a particular type of partition having block
size at most two and exactly two, respectively. Their standard
representation is a unique graph on the vertex set $[N]$ whose edge
set consists of arcs connecting the elements of each block in
numerical order, see Fig.~\ref{F:parcro1}.

Given a (set) partition $P$, a \emph{$k$-crossing} is a set of $k$
edges $\{(i_{1},j_{1}),(i_{2},j_{2}),\ldots, (i_{k},j_{k})\}$ such
that $i_{1}< i_{2}< \ldots < i_{k}< j_{1}< j_{2}< \ldots < j_{k}$,
see Fig.\ref{F:parcro1}. A \emph{$k$-noncrossing partition} is a
partition without any $k$-crossings.  We denote the sets of
$k$-noncrossing partitions and $m$-regular, $k$-noncrossing
partitions by $\mathcal{P}_k(N)$ and $\mathcal{P}_{k,m}(N)$,
respectively. For instance, the set of $2$-regular, $3$-noncrossing
partitions are denoted by $\mathcal{P}_{3,2}(N)$.

\begin{figure}[ht]
\centerline{\includegraphics[width=0.75\textwidth]{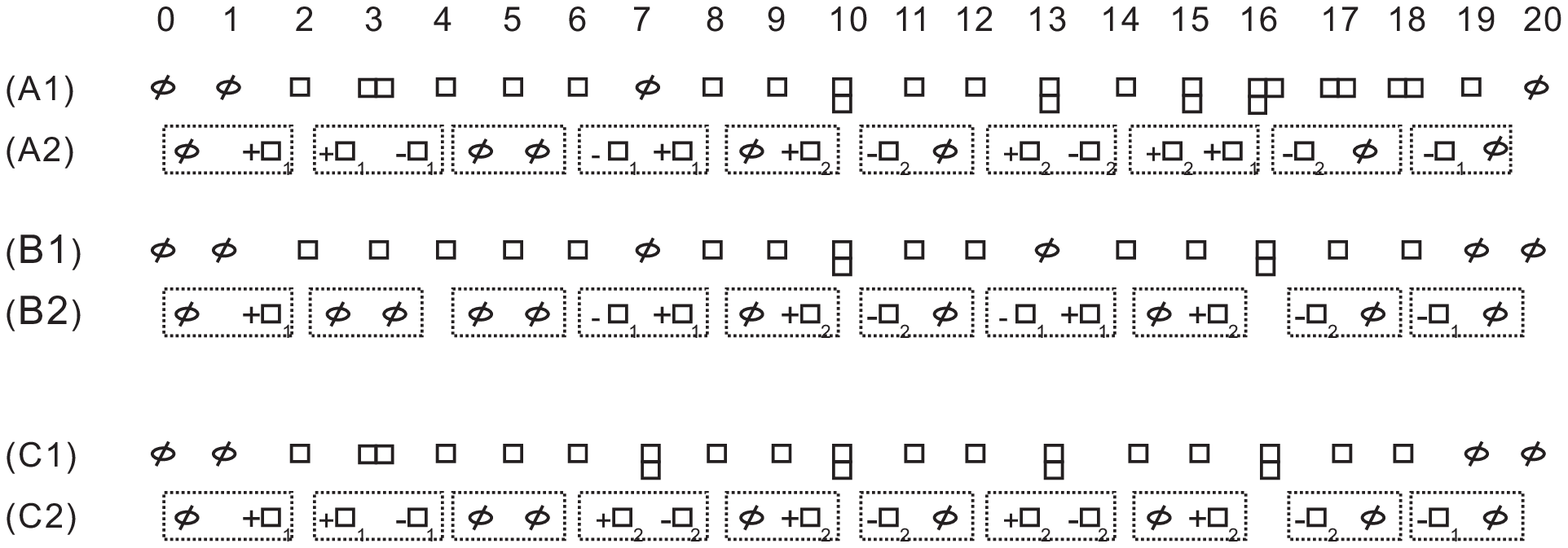}}
\caption{\small
 Vacillating tableaux and elementary moves:
\textsf{(A1)} shows a general vacillating tableaux with EM
\textsf{(A2)}. \textsf{(B1)} displays another vacillating tableaux
and its set of EM, $\{(-\square,\varnothing),
(\varnothing,+\square),(\varnothing,\varnothing),(-\square,+\square)\}$
shown in \textsf{(B2)}. In \textsf{(C1)} we present a vacillating
tableaux with the EM $\{(-\square,\varnothing),
(\varnothing,+\square),(\varnothing,\varnothing),(+\square,-\square)\}$
displayed in \textsf(C2).} \label{F:tableaux}
\end{figure}

A \emph{(generalized)} vacillating tableau \cite{Reidys:08tan}
$V_{\lambda}^{2N}$ of shape $\lambda$ and length $2N$ is a sequence
$\lambda^{0}, \lambda^{1},\ldots,\lambda^{2N}$ of shapes such that
{\sf (1)} $\lambda^{0}=\varnothing$, $\lambda^{2N}=\lambda$ and {\sf
(2)} for $1\le i\le N$, $\lambda^{2i-1},\lambda^{2i}$ are derived
from $\lambda^{2i-2}$ by \emph{elementary moves} (EM) defined as
follows: $(\varnothing,\varnothing)$: do nothing twice;
$(-\square,\varnothing)$: first remove a square then do nothing;
$(\varnothing,+\square)$: first do nothing then add a square; $(\pm
\square,\pm \square)$: add/remove a square at the odd and even
steps, respectively. We use the following notation: if
$\lambda_{i+1}$ is obtained from $\lambda_{i}$ by adding, removing a
square from the $j$-th row, or doing nothing we write
$\lambda_{i+1}\setminus\lambda_{i}=+\square_{j}$,
$\lambda_{i+1}\setminus\lambda_{i}= -\square_{j}$ or
$\lambda_{i+1}\setminus\lambda_{i}=\varnothing$, respectively, see
Fig.\ref{F:tableaux}.

\begin{figure}[ht]
\centerline{%
\epsfig{file=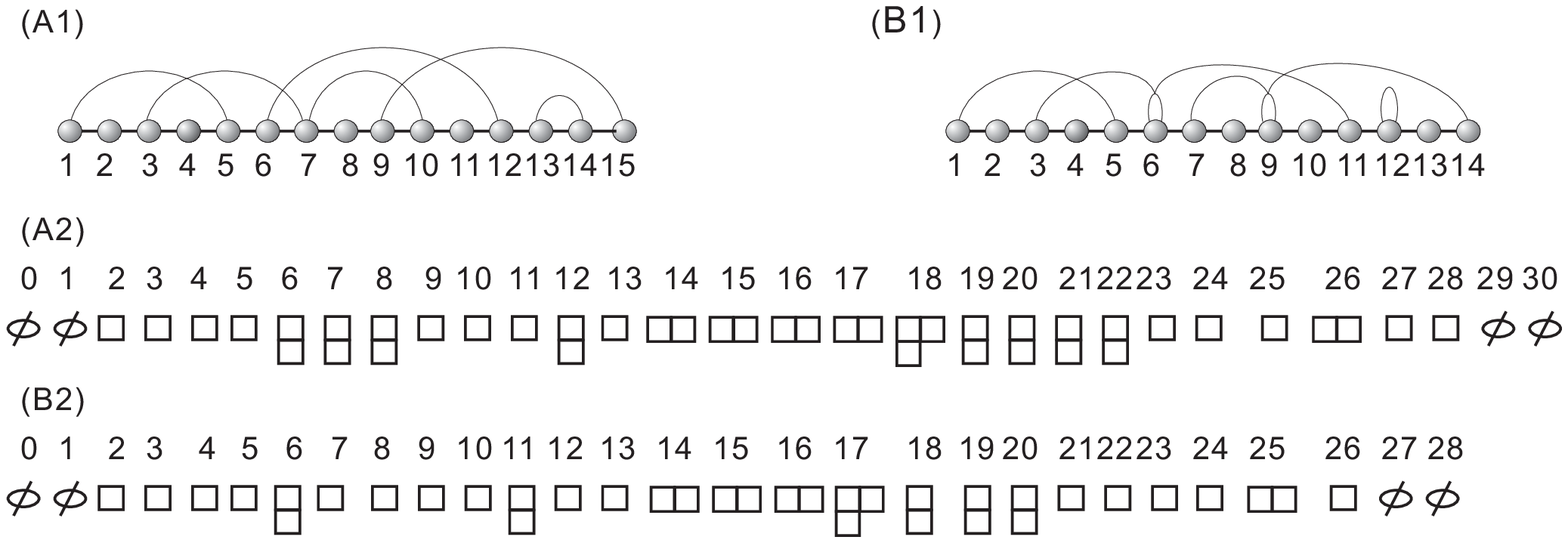,width=0.9\textwidth}\hskip15pt
 }
\caption{\small Vacillating tableaux, partitions and braids: in
\textsf{(A1)} we show a $3$-noncrossing partition and in
\textsf{(A2)} the associated vacillating tableau. \textsf{(B1)}
shows a $3$-noncrossing braid and  \textsf{(B2)} its vacillating
tableau.} \label{F:braid}
\end{figure}

A braid over $[N]$ can be represented via introducing loops $(i,i)$
and drawing arcs $(i,j)$ and $(j,\ell)$ with $i<j<\ell$ as crossing,
see Fig.\ref{F:braid} \textsf{(B1)}. A \emph{$k$-noncrossing braid}
is a braid without any $k$-crossings. We denote the set of
$k$-noncrossing braids over $[N]$ with and without isolated points
by $\mathcal{B}_k(N)$ and $\mathcal{B}^{*}_k(N)$, respectively. Chen
{\it et al.} \cite{Chen} have shown that each $k$-noncrossing
partition corresponds uniquely to a vacillating tableau of empty
shape, having at most $(k-1)$ rows, obtained via the EM
$\{(-\square, \varnothing),(\varnothing, +\square), (\varnothing,
\varnothing),(-\square, +\square)\}$, see Fig.\ref{F:braid}
\textsf{(A2)}. In \cite{Reidys:08tan}, Chen {\it et al.} proceed by
proving that vacillating tableaux of empty shape, having at most
$(k-1)$-rows which are obtained by the EMs $\{(-\square,
\varnothing),(\varnothing, +\square),(\varnothing, \varnothing),
(+\square, -\square)\}$ correspond uniquely to $k$-noncrossing
\emph{braids}.

In the following, we consider $3$-noncrossing partition or a
$2$-regular, $3$-noncrossing partition of $N$ vertices, where $N
\geq 1$. To this end, we introduce two $\mathbb{Z}^2$ domains
$Q_2=\{(a_1,a_2)\in\mathbb{Z}_2\mid a_1, a_2\geq 0\}$ and $W_2=
\{(a_1,a_2)\in\mathbb{Z}_2\mid a_1>a_2\geq 0\}$. Let $D\in
\{Q_2,W_2\}$, $\mathbf{e_1} = (1,0)$ and $\mathbf{e_{2}}=(0,1)$. A
\emph{$\mathcal{P}_D$-walk} is a lattice walk in $D$ starting at
$(1,0)$ having steps $\pm \mathbf{e_{1}}$, $\pm \mathbf{e_{2}}$,
$\mathbf{0}=(0,0)$, such that even steps are $+\mathbf{e_{1}}$,
$+\mathbf{e_{2}}$, or $\mathbf{0}$ and odd steps are either
$-\mathbf{e_{1}}$, $-\mathbf{e_{2}}$, or $\mathbf{0}$.
Analogously, a \emph{$\mathcal{B}_D$-walk} is a lattice walk in $D$
starting at $(1,0)$ whose even steps are either $-\mathbf{e_{1}}$,
$-\mathbf{e_{2}}$, or $\mathbf{0}$ and whose odd steps are
$+\mathbf{e_{1}}$, $+\mathbf{e_{2}}$, or $\mathbf{0}$. By abuse of
language, we will omit the subscript $D$.

Interpreting the number of squares in the rows of the shapes
as coordinates of lattice points, we immediately obtain
\begin{theorem}\label{T:bijection}\cite{Reidys:08tan}\\
{\rm (1)} The number of 3-noncrossing partitions over $[N]$ equals
the number of $\mathcal{P}_{W_2}$-walks from $(1,0)$ to
itself of length $2N$. \\
{\rm (2)} The number of 3-noncrossing braids over $[N]$ equals the
number of $\mathcal{B}_{W_2}$-walks from $(1,0)$ to itself of length
$2N$.
\end{theorem}

\section{Random $3$-noncrossing partitions}\label{S:vacillating}

In the following, set $\bar{x}=1/x$ and $\bar{y}=1/y$. We work in
the ring of power series in $t$ whose coefficients are Laurent
polynomials in $x$ and $y$, denoted by
$Q[x,\bar{x},y,\bar{y}][[t]]$. We first review some calculations of
\cite{Xin:06}. Let $a^{i,j}_{s}$ denote the number of
$\mathcal{P}_{Q_2}$-walks of length $s$ ending at $(i,j)\in Q_2$ and
\begin{equation}\label{E:XY}
{F}_{e}(x,y;t)=\sum_{i,j,\ell}a^{i,j}_{2\ell}x^iy^jt^{2\ell}  \quad
\text{\rm and}\quad
{F}_{o}(x,y;t)=\sum_{i,j,\ell}a^{i,j}_{2\ell+1}x^iy^jt^{2\ell+1}
\end{equation}
denote the generating functions of $\mathcal{P}_{Q_{2}}$-walks of
even and odd length and ending at $(i,j)$, respectively. In
particular, let $H_{e}(x; t)=F_{e}(x,0;t)$ and
$V_{e}(y;t)=F_{e}(0,y;t)$ denote the generating functions of even
$\mathcal{P}_{Q_{2}}$-walks ending on the $x$-axis and $y$-axis.
\begin{proposition}\label{P:kernel1}\cite{Xin:06}
Let $V(y; t)= V_e(y;t^2)$, $H(x;t)= H_e(x;t^2)$, then
\begin{eqnarray}\label{E:define}
&&F_{e}(x,y;t)=x+F_{o}(x,y;t)(1+x+y)t,\\
\label{E:kernel}
&&K(x,y;t)\, F(x,y;t)=xy+x^2y+x^2-xH(x,t)-yV(y,t),
\end{eqnarray}
where $F(x,y;t)=\sum_{i,j,\ell} f^{i,j}_{\ell}x^{i}y^{j}t^{\ell}$ is
given by $F_o(x,y;t)= tF(x,y;t^2)$ and
\begin{equation}
K(x,y;t):=xy-t(1+x+y)(x+y+xy)
\end{equation}
is called the kernel of eq.~(\ref{E:kernel}).
\end{proposition}

Eq.~(\ref{E:define}) follows from $a^{1,0}_{0}=1$ and
$a^{i,j}_{2\ell}=a^{i-1,j}_{2\ell-1}+
a^{i,j-1}_{2\ell-1}+a^{i,j}_{2\ell-1}$ for $1 \leq \ell\leq  N$. In
view of Proposition \ref{P:kernel1}, $F(x,y;t)$ is the key for
enumerating $\mathcal{P}_{Q_2}$-walks. We observe that
\begin{equation*}
Y_{0}=\frac{1-(\bar{x}+3+x)t-\sqrt{(1-(1+x+\bar{x}t))^2-4t}}
{2(1+\bar{x})t}
\end{equation*}
is the unique power-series solution in $t$ of the
kernel equation $K(x,y;t)=0$, whence
\begin{equation}
xY_{0}=t(1+x+Y_{0})(x+Y_{0}+xY_{0}).
\end{equation}
Since $K(x,y;t)$ is quadratic in both $x$ and $y$, $(x,Y_{0})$ can
be regarded as a pair of Laurent series in $t$ with coefficients
satisfying $K(x,Y_{0};t)= 0$. As $(\bar{x}Y_{0},Y_{0})$ and
$(\bar{x}Y_{0},\bar{x})$ also solve $K(x,y;t)=0$, we arrive at
\begin{eqnarray}\label{E:1}
xH(x;t) + Y_{0}V(Y_{0};t) &=& xY_{0} + x^2Y_{0} + x^2,\\\label{E:2}
\bar{x}Y_{0}H(\bar{x}Y_{0};t)+Y_{0}V(Y_{0};t)&=&\bar{x}Y_{0}^2 +
\bar{x}^2Y_{0}^3 +\bar{x}^2Y_{0}^2,\\\label{E:3}
\bar{x}Y_{0}H(\bar{x}Y_{0};t)+\bar{x}V(\bar{x};t)&=&\bar{x}^2Y_{0}+
\bar{x}^3Y_{0}^2 +\bar{x}^2Y_{0}^2.
\end{eqnarray}
Indeed, we can conclude from eqs.~(\ref{E:1})-(\ref{E:3})
\begin{equation}\label{E:obj}
xH(x)+ \bar{x}V (\bar{x}) = x^2 + (\bar{x}^2 + x + x^2)Y_{0} +
(\bar{x}^3-\bar{x})Y_{0}^2-\bar{x}^2Y_{0}^3.
\end{equation}
Since $H$ and $V$ are generating functions, $xH(x;t)$ contains only
positive powers of $x$. Similarly, $\bar{x}V(\bar{x};t)$ contains
only negative powers of $x$. Therefore, setting ${\sf NT}_{x}$ and
${\sf PT}_{x}$ to be the operators that extract the positive and
negative powers of $x$ from a power series in
$\mathbb{Q}[x,\bar{x}][[t]]$, we arrive at
\begin{eqnarray}
\label{E:vh1}xH(x;t)&=&{\sf
PT}_{x}[x^2+(\bar{x}^{2}+x+x^2)Y_{0}+(\bar{x}^{3}-\bar{x})
Y_{0}^{2}-\bar{x}^{2}Y_{0}^{3}]\\
\label{E:vh2} \bar{x}V(\bar{x};t)&=&{\sf
NT}_{x}[x^2+(\bar{x}^{2}+x+x^2)Y_{0}+(\bar{x}^{3}-\bar{x})
Y_{0}^{2}-\bar{x}^{2}Y_{0}^{3}].
\end{eqnarray}
Via Lagrange inversion \cite{ec1}, we derive the formula
\begin{equation}\label{E:inverse1}
[x^{0}t^s]x^{\ell}Y_{0}^{m}=\sum_{j}\frac{m}{s}{s \choose j} {s
\choose {j+m}}{{2j+m}\choose {j-\ell}}.
\end{equation}
The enumeration of $\mathcal{P}_{W_2}$-walks follows from the
reflection principle \cite{Gessel:92a}. Let $\omega_{\ell}^{i,j}$
denote the number of $\mathcal{P}_{W_2}$-walks ending at $(i,j)$ of
length $\ell$. The reflection-principle implies
\begin{equation}\label{E:reflection1}
\omega^{i,j}_{\ell}=a^{i,j}_{\ell}-a^{j,i}_{\ell}.
\end{equation}
In the following lemma, we express $\omega_{\ell}^{i,j}$ via the
coefficients of $F(x,y;t)$ and establish their corresponding
recurrence relations.
\begin{lemma}\label{L:p2}
{\rm (a)} Suppose $1\leq \ell\leq N-1$ and $(i,j)\in\mathbb{Z}^2$,
then $\omega_{s}^{i,j}$ is given by
\begin{equation}\label{E:look0}
\omega^{i,j}_{s}=
\begin{cases}
f^{i,j}_{\ell}-f^{j,i}_{\ell} & \text{\rm for $s=2\ell+1$},\\
f^{i,j}_{\ell}+f^{i-1,j}_{\ell}+f^{i,j-1}_{\ell}- f^{j,i}_{\ell}-
f^{j-1,i}_{\ell}-
f^{j,i-1}_{\ell} &  \text{\rm for $s=2\ell+2$},\\
\end{cases}
\end{equation}
where {\sf (1)} $\omega^{i,j}_{0}=1$, for $i=1,j=0$ and
$\omega^{i,j}_{0}=0$, otherwise; {\sf (2)} $f^{i,j}_{\ell}=0$ for
$i\not\in \{0,1,\dots, \ell\}$ or $j\not\in \{0,1,\dots, \ell\}$. \\
{\rm (b)} $f^{i,j}_{\ell}$ satisfies the recursion
\begin{equation}\label{E:hh}
f^{i,j}_{\ell}=
\begin{cases}
[y^{j}t^{\ell+1}]V(y;t)-f^{0,j-1}_{\ell}=
\sum_{j_1}(-1)^{j-j_1}[y^{j_1}t^{\ell+1}]V(y;t)
&  \text{\rm for $i=0$, $j\neq 0$}\\
[x^{i}t^{\ell+1}]H(x;t)-f^{i-1,0}_{\ell}=
\sum_{i_1}(-1)^{i-i_1}[x^{i_1}t^{\ell+1}]H(x;t)
& \text{\rm for $i\neq 0$, $j=0$}\\
f^{i,j+1}_{\ell-1}-f^{i+1,j}_{\ell-1}-f^{i-1,j+1}
_{\ell-1}-f^{i+1,j-1}_{\ell-1}-
3f^{i,j}_{\ell-1}-f^{i-1,j}_{\ell-1}-f^{i,j-1}_{\ell-1} & \text{\rm
otherwise.}
\end{cases}
\end{equation}
\end{lemma}
\begin{proof}
We first prove assertion (a). Indeed, according to
Proposition~\ref{P:kernel1} we have $F_{o}(x,y;t)=tF(x,y;t^2)$ and
$F_{e}(x,y;t)= x+F_{o}(x,y;t)(1+x+y)t$ or equivalently
\begin{eqnarray}
a^{i,j}_{2\ell+1}= f^{i,j}_{\ell}  \quad \text{\rm and} \quad
a^{i,j}_{2\ell}
=f^{i,j}_{\ell-1}+f^{i-1,j}_{\ell-1}+f^{i,j-1}_{\ell-1},
\end{eqnarray}
whence eq.~(\ref{E:look0}) follows from
$\omega^{i,j}_{\ell}=a^{i,j}_{\ell}-a^{j,i}_{\ell}$.\\
Next, we prove (b). We distinguish the following three cases. First,
suppose $i\neq 0$ and $j\neq 0$. Equating the coefficients of
$x^iy^jt^{\ell}$ in eq.~(\ref{E:kernel}), we derive
\begin{equation}
f^{i,j}_{\ell}=
f^{i,j+1}_{\ell-1}-f^{i+1,j}_{\ell-1}-f^{i-1,j+1}_{\ell-1}-f^{i+1,j-1}_{\ell-1}-
3f^{i,j}_{\ell-1}-f^{i-1,j}_{\ell-1}-f^{i,j-1}_{\ell-1}.
\end{equation}
In case of $i=0$ or $j=0$, let $[x^it^{\ell}]xH(x;t)$ and
$[y^jt^{\ell}]yV(y;t)$ denote the coefficients of $x^it^{\ell}$ in
$xH(x;t)$ and the coefficient of $y^{j}t^{\ell}$ in $yV(y;t)$,
respectively. Then we have
\begin{eqnarray}\label{E:V}
f_{{\ell}-1}^{0,j-1}+f_{{\ell}-1}^{0,j-2}&=&[y^jt^{\ell}]yV(y;t),\\
\label{E:H}
f_{{\ell}-1}^{i-1,0}+f_{{\ell}-1}^{i-2,0}&=&[x^it^{\ell}]xH(x;t).
\end{eqnarray}
According to eq.~(\ref{E:vh1}) and eq.~(\ref{E:vh2}), the
coefficients $[x^it^{\ell}]xH(x;t)$ and $[y^jt^{\ell}]yV(y;t)$ are
given as follows:
\begin{eqnarray}\label{E:xh}
[x^{i}t^{\ell}]H(x,t)&=&[x^{0}t^{\ell}]{x}^{-i-1}(x^2+
(\bar{x}^{2}+x+x^2)
Y_{0}+(\bar{x}^{3}-\bar{x})Y_{0}^{2}-\bar{x}^{2}Y_{0}^{3}),\\
\label{E:xh2} [y^{j}t^{\ell}]V(y,t)&=&[x^{0}t^{\ell}]{x}^{j+1}(x^2+
(\bar{x}^{2}+x+x^2)Y_{0}+(\bar{x}^{3}-\bar{x})Y_{0}^{2}-\bar{x}^{2}Y_{0}^{3}).
\end{eqnarray}
whence (b) and the proof of the lemma is complete.
\end{proof}
Lemma \ref{L:p2} allows us to compute $\omega^{i,j}_{\ell}$ for all
$i,j\in \mathbb{Z}^2$ and $1\leq \ell\leq 2N$. In the following, we
consider a vacillating tableaux as the sampling path of a
Markov-process, whose transition probabilities can be calculated via
the terms $\omega^{i,j}_{\ell}$.

\begin{theorem}\label{T:algpar}
Algorithm~$1$ generates a random $3$-noncrossing partition after a
pre-processing step having $O(N^4)$ time and $O(N^3)$ space
complexity with uniform probability in linear time and space.
\end{theorem}
\begin{algorithm}\label{A:algor1}
\begin{algorithmic}[1]
\STATE{$i$=1}

\STATE{\it Tableaux} (Initialize the sequence of shapes,
$\{\lambda^{i}\}_{i=0}^{i=2N}$)

\STATE {$\lambda^{0}=\varnothing$, $\lambda^{2N}=\varnothing$}

\WHILE {$i < 2N$ }

\IF {$i$ is even}

\STATE  X[0]$\leftarrow$ ${\sf V}_3(\lambda^{i+1}_{\varnothing},
2N-(i+1))$

\STATE  X[1]$\leftarrow$ ${\sf V}_3(\lambda^{i+1}_{-\square_{1}},
2N-(i+1))$

\STATE  X[2]$\leftarrow$ ${\sf V}_3(\lambda^{i+1}_{-\square_{2}},
2N-(i+1))$

\ENDIF

\IF {$i$ is odd}

\STATE  X[0]$\leftarrow$ ${\sf V}_3(\lambda^{i+1}_{\varnothing},
2N-(i+1))$

\STATE  X[1]$\leftarrow$ ${\sf V}_3(\lambda^{i+1}_{+\square_{1}},
2N-(i+1))$

\STATE  X[2]$\leftarrow$ ${\sf V}_3(\lambda^{i+1}_{+\square_{2}},
2N-(i+1))$

\ENDIF

\STATE {\it sum} $\leftarrow$ X[0]+X[1]+X[2]

\STATE {\it Shape} $\leftarrow$ {\sf Random}({\it sum}) ({\sf
Random} generates the random shape $\lambda^{i+1}_{+\square_{j}(or
-\square_{j})}$ with probability ${X[j]}/{\it sum}$ or
$\lambda^{i+1}_{\varnothing}$ with probability $X[0]/sum$)

\STATE $i \leftarrow i+1$

\STATE Insert {\it Shape} into {\it Tableau} (the sequence of
shapes).

\ENDWHILE

\STATE  {\sf Map}({\it Tableau}) (maps {\it Tableau} into its
corresponding $3$-noncrossing partition)
\end{algorithmic}
\caption{\small Uniform generation of $3$-noncrossing partitions}
\end{algorithm}
\begin{proof}
The main idea is to interpret tableaux of $3$-noncrossing partitions
as sampling paths of a stochastic process. We label the $(i+1)$-th
shape, $\lambda^{i+1}_{\alpha}$ by $\alpha=\lambda^{i+1}\setminus
\lambda^{i}\in \{+\square_{1}, +\square_{2}, -\square_{1},
-\square_{2}, \varnothing\}$, where the labeling specifies the
transition from $\lambda^i$ to
$\lambda^{i+1}$.\\
Let ${\sf V}(\lambda^{i+1}_{\alpha},2N-(i+1))$ denote the number of
vacillating tableaux of length $(i+1)$ such that
$\lambda^{i+1}\setminus \lambda^{i}=\alpha$. We remark here that if
$\lambda^{i+1}$ has \NEW{$(a-1)$} in the first row and $b$ boxes in
the second row, then ${\sf V}(\lambda^{i+1}_{\alpha},2N-(i+1))=
\omega^{a,b}_{2N-(i+1)}$.
\begin{figure}[ht]
\centerline{\includegraphics[width=0.6
\textwidth]{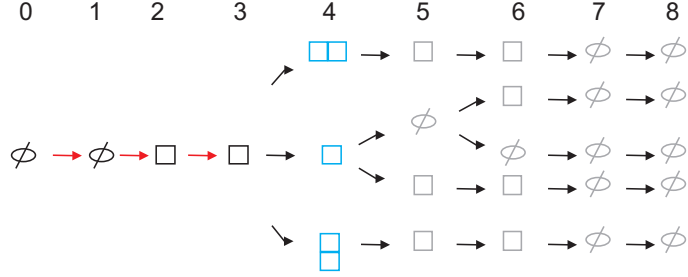}} \caption{\small Sampling paths: we
display all tableaux sequences where $\lambda_1=\varnothing$,
$\lambda^2=\lambda^3=\square$. We display the three possible
$\lambda^4$-shapes (blue), induced by $\varnothing,+
\square_1,+\square_2$ and all possible continuations (grey).}
\label{F:rest}
\end{figure}
Let $(X^i)_{i=0}^{2N}$ be given as follows:\\
$\bullet$ $X^0=X^{2N}=\varnothing$ and $X^i$ is a shape having at
most 2 rows;\\
$\bullet$ for $1\le i\le N-1$, we have $X_{2i+1}\setminus X_{2i}\in
\{\varnothing, -\square_{1}, -\square_{2}\}$ and $X_{2i+2}\setminus
X_{2i+1}\in
\{\varnothing, +\square_{1}, +\square_{2}\}$.\\
$\bullet$ for $1\leq i\leq 2N-1$, we have
\begin{equation}\label{E:w0}
\mathbb{P}(X^{i+1}=\lambda^{i+1}\mid X^i=\lambda^i)= \frac{{\sf
V}(\lambda^{i+1}_{\alpha},2N-i-1)}{{\sf V}(\lambda^{i},2N-i)}.
\end{equation}
In view of eq.~(\ref{E:w0}), we immediately observe
\begin{equation}
\prod_{i=0}^{2N-1}\mathbb{P}(X^{i+1}=\lambda^{i+1}\mid
X^i=\lambda^i)= \frac{{\sf V}(\lambda^{2N}=\varnothing,0)} {{\sf
V}(\lambda^{0}=\varnothing,N)} =\frac{1}{{\sf V}(\varnothing,2N)}.
\end{equation}
Consequently, the process $(X^i)_{i=0}^{2N}$ generates random
$3$-noncrossing partitions with uniform probability in $O(N)$ time and space.\\
As for the derivation of the transition probabilities, suppose we
are given a shape $\lambda^{h}$ having exactly $a$ and $b$ squares
in the first and second row, respectively. According to
Lemma~\ref{L:p2}, we obtain
\begin{equation*}
\omega^{a,b}_{2N-h}=
\begin{cases}
f^{a,b}_{\ell}-f^{b,a}_{\ell} &  \text{\rm for $2N-h=2\ell+1$},\\
f^{a,b}_{\ell}+f^{a-1,b}_{\ell}+f^{a,b-1}_{\ell}- f^{b,a}_{\ell}-
f^{b-1,a}_{\ell}- f^{b,a-1}_{\ell} &\text{\rm for $2N-h=2\ell+2$},
\end{cases}
\end{equation*}
which reduces the problem to the calculation of at most six
coefficients $f^{i,j}_{\ell}$ for fixed $i$, $j$ and $\ell$. \\
{\it Claim.} The coefficients $f^{i,j}_{\ell}$ for all $0\le i,j\leq
\ell$ and $1\leq \ell\leq N-1$,
can be computed in $O(N^4)$ time and $O(N^3)$ space complexity.\\
{\sf Step 1.} For fixed $i$, $j$ and $\ell$, $f^{i,0}_{\ell}$ and
$f^{0,j}_{\ell}$ can be computed in $O(N^2)$ time. Indeed, according
to eq.~(\ref{E:hh}), we have
\begin{equation*}
f^{i,0}_{\ell}=\sum_{i_1}(-1)^{i-i_1}[x^{i_1}t^{\ell+1}]H(x;t),
\end{equation*}
whence $f^{i,0}_{\ell}$ can be calculated via a {\bf For}-loop
summing over the terms $(-1)^{i-i_1}[x^{i_1}t^{\ell+1}]H(x;t)$. In
view of eq.~(\ref{E:xh2}), the $[x^{i_1}t^{\ell+1}]H(x;t)$ terms in
turn are expressed via $[x^{0}t^{\ell+1}]x^{q}Y_{0}^{m}$, for $q\in
\{-i_1-4,-i_1-3,-i_1-2,-i_1,-i_1+1\}$ and $m\in \{0,1,2\}$.
According to eq.~(\ref{E:inverse1}), the quantities
$[x^{0}t^{\ell+1}]x^{q} Y_{0}^{m}$ can be calculated via a {\bf
For}-loop summing over the terms $\frac{m}{{\ell+1}}{{\ell+1}
\choose j}{{\ell+1} \choose {j+m}}{{2j+m}\choose {j-q}}$ for $1\leq
j\leq {\ell+1}$. Consequently, for fixed $i,\ell$ the coefficient
$f^{i,0}_{\ell}$ can be derived in $O(N^2)$ time. Using the same
arguments we obtain for fixed
$j,\ell$ the coefficient $f^{0,j}_{\ell}$ in $O(N^2)$ time.\\
{\sf Step 2.} We compute  $f^{i,0}_{\ell}$ and $f^{0,j}_{\ell}$ for
all $i,j$
and $\ell$ via two nested {\bf For}-loops in $O(N^2)$ time.\\
{\sf Step 3.} Once the coefficients $f^{i,0}_{\ell}$ and
$f^{0,j}_{\ell}$ are calculated for all $i,j,\ell$, we compute
$f^{i,j}_{\ell}$ for arbitrary $i,j,\ell$ employing three nested
{\bf For}-loops and the recurrence of eq.~(\ref{E:hh}) since there
exists boundary conditions as follows:\\
\begin{eqnarray}
f^{i,j}_{\ell}&=&0\quad i>\ell\ \text{or}\ j>\ell,\\
f^{i,j}_{0}&=&
\begin{cases}
1\quad i=1,\ j=0;\\
0\quad \text{otherwise}.
\end{cases}
\end{eqnarray}
Therefore, we compute $f^{i,j}_{\ell}$ for arbitrary $i,j,\ell$ in
$O(N^4)$ time and $O(N^3)$ space and the claim follows.\\
Given the coefficients $f^{i,j}_{\ell}$ for all $i,j,\ell$, we can
derive the transition probabilities in $O(1)$ time. Accordingly, we
obtain all transition probabilities ${\sf
V}(\lambda^i_{\alpha},2N-i))$ in $O(N^4)$ time and $O(N^{3})$ space
complexity.
\end{proof}
\section{Random $2$-regular, $3$-noncrossing partitions}\label{S:braid}

In this section, we generate random $2$-regular, $3$-noncrossing
partitions employing a bijection between the set of $2$-regular,
$3$-noncrossing partitions over $[N]$ and the set of $k$-noncrossing
braids without loops over $[N-1]$.

\begin{lemma}\label{L:bijection-2}
Set $k\in\mathbb{N}$, $k\ge 3$. Suppose
$(\lambda^i)_{i=0}^{2\ell+1}$ in which $\lambda^{0}=\varnothing$ is
a sequence of shapes such that
$\lambda^{2j}\setminus\lambda^{2j-1}\in
\{\{+\square_{h}\}_{h=1}^{h=k-1},\varnothing\}$ and
$\lambda^{2j-1}\setminus \lambda^{2j-2}\in
\{\{-\square_{h}\}_{h=1}^{h=k-1},\varnothing\}$. Then
$(\lambda^i)_{i=0}^{2\ell+1}$ induces a unique sequence of shapes
$(\mu^i)_{i=0}^{2\ell}$ with the following properties
\begin{eqnarray}
&& \label{E:right}
 \mu^{2j+1}\setminus
\mu^{2j}\in \{\{+\square_{h}\}_{h=1}^{h=k-1},\varnothing\}\
\text{\it and}\
\mu^{2j+2}\setminus\mu^{2j+1}\in \{\{-\square_{h}\}_{h=1}
^{h=k-1},\varnothing\},\\
\label{E:shape} &&
 \mu^{2j}=\lambda^{2j+1}\\\label{E:shape1}
&&  \mu^{2j+1}\neq\lambda^{2j+2}\quad \Longrightarrow\quad
\mu^{2j+1}\in \{\lambda^{2j+1},\lambda^{2j+3}\} .
\end{eqnarray}
\end{lemma}
\begin{figure}[ht]
\centerline{%
\epsfig{file=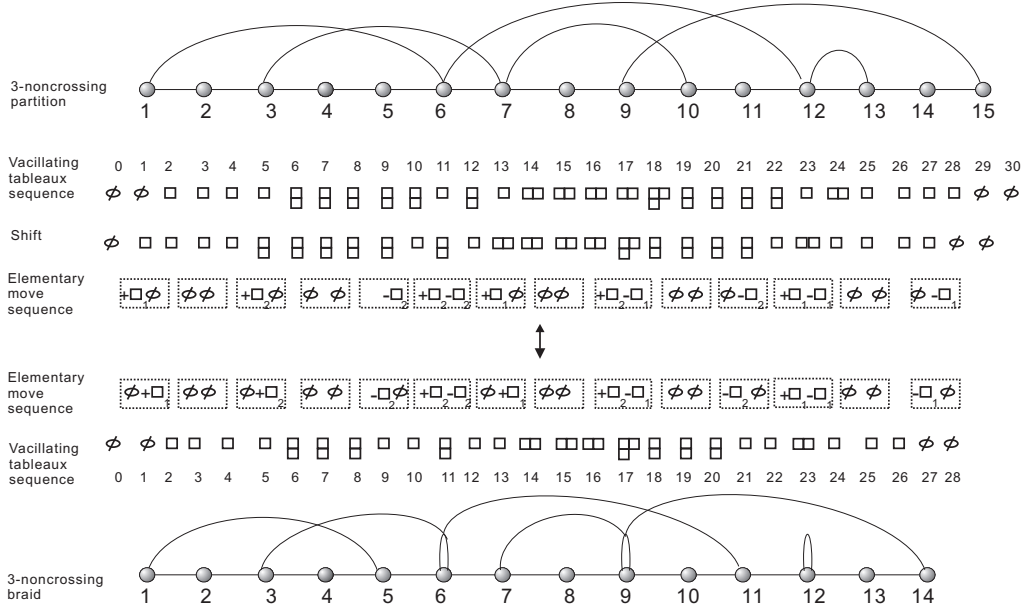,width=0.9\textwidth}\hskip15pt
 }
\caption{\small Illustration of Lemma~\ref{L:bijection-2}. We show
how to map the $\mathcal{P}_{W_2}$-walk induced by the
$3$-noncrossing partition $P$, of Fig.~\ref{F:braid}, \textsf{(A)}
into the $\mathcal{B}_{W_2}$-walk. The latter corresponds to the
$3$-noncrossing braid of Fig.~\ref{F:braid}, \textsf{(B)}. Here each
vertex, $i$, is aligned with the triple of shapes
$(\lambda^{2i-2},\lambda^{2i-1},\lambda^{2i})$ in the corresponding
tableaux.} \label{F:map}
\end{figure}
\begin{proof}
Since $\lambda^1=\varnothing$, $(\lambda^j)_{j=1}^{2\ell+1}$ corresponds
to a sequence of pairs $((x_i,y_i))_{i=1}^{\ell}$ given by
$x_{i}=\lambda^{2i}\setminus\lambda^{2i-1}$ and
$y_{i}=\lambda^{2i+1}\setminus\lambda^{2i}$ such that
\begin{equation}\label{E:neu}
\forall\, 1\le i\le \ell;\qquad (x_i,y_i)\in
\{(\varnothing,\varnothing),\  (+\square_h,\varnothing), \
(\varnothing,-\square_h), \ (+\square_h,-\square_j)\}.
\end{equation}
Let $\varphi$ be given by
\begin{equation}\label{E:varphi2}
\varphi(({x}_i,{y}_i))=
\begin{cases}
({x}_i,{y}_i) & \ \text{\rm for } \
({x}_i,{y}_i)=(+\square_h,-\square_j) \\
({y}_i,{x}_i) & \ \text{\rm otherwise,}
\end{cases}
\end{equation}
and set $(\varphi(x_i,y_i))_{i=1}^\ell=(a_i,b_i)_{i=1}^{\ell}$. Note
that $(a_i,b_i)\in
\{(-\square_h,\varnothing),(\varnothing,+\square_h),
(\varnothing,\varnothing),(+\square_h,-\square_j)\}$, where $1\leq
h,j\leq k-1$. Let $(\mu^i)_{i=1}^{\ell}$ be the sequence of shapes
induced by $(a_i,b_i)_{i=1}^{\ell}$ according to
$a_{i}=\mu^{2i}\setminus\mu^{2i-1}$ and
$b_{i}=\mu^{2i-1}\setminus\mu^{2i-2}$ initialized with
$\mu_{0}=\varnothing$. Now eq.~(\ref{E:right}) is implied by
eq.~(\ref{E:neu}) and eq.~(\ref{E:shape}) follows by construction.
Suppose $\mu^{2j+1}\neq\lambda^{2j+2}$ for some $0\le j\le \ell-1$.
By definition of $\varphi$, only pairs containing ``$\varnothing$''
in at least one coordinate are transposed from which we can conclude
$\mu^{2j+1}=\mu^{2j}$ or $\mu^{2j+1}=\mu^{2j+2}$, whence
eq.~(\ref{E:shape1}). I.e.~we have the following situation
\begin{equation*}
\diagram
& \lambda^{2j+1} \ar@{=}[dl]\rto^{} & \lambda^{2j+2}\ar@{-}[dl]
\rto^{} &
\lambda^{2j+3} \ar@{=}[dl]  \\
 \mu^{2j} \ar@{=}[r] & \mu^{2j+1} \rto^{} &\mu^{2j+2}
\enddiagram
 \text{\rm or }
\diagram
& \lambda^{2j+1} \ar@{=}[dl]\rto^{} & \lambda^{2j+2}\ar@{-}[dl]
\rto^{} &
\lambda^{2j+3} \ar@{=}[dl]  \\
 \mu^{2j} \rto^{} & \mu^{2j+1} \ar@{=}[r] &\mu^{2j+2},
\enddiagram
\end{equation*}
and the lemma follows.
\end{proof}
Lemma~\ref{L:bijection-2} establishes a bijection between
$\mathcal{P}_{W_{k-1}}$-walks of length $2\ell+1$ and
$\mathcal{B}_{W_{k-1}}$-walks of length $2\ell$, where
$W_{k-1}=\{(a_1,a_2,\dots,a_{k-1})\mid a_1>a_2>\dots>a_{k-1}\}$.
\begin{corollary}\label{C:bi}
Let $\mathcal{P}_k(N)$ and $\mathcal{B}_{k}(N-1)$ denote the set of
$k$-noncrossing partitions on $[N]$ and $k$-noncrossing braids on
$[N-1]$.
Then \\
{\rm (a)} we have a bijection
\begin{equation}\label{E:biject2}
\vartheta\colon \mathcal{P}_{k}(N)\longrightarrow
\mathcal{B}_k(N-1).
\end{equation}
{\rm (b)} $\vartheta$ induces by restriction a bijection between
$2$-regular, $k$-noncrossing partitions on $[N]$ and $k$-noncrossing
braids without loops on $[N-1]$.
\end{corollary}

\begin{figure}[ht]
\centerline{%
\epsfig{file=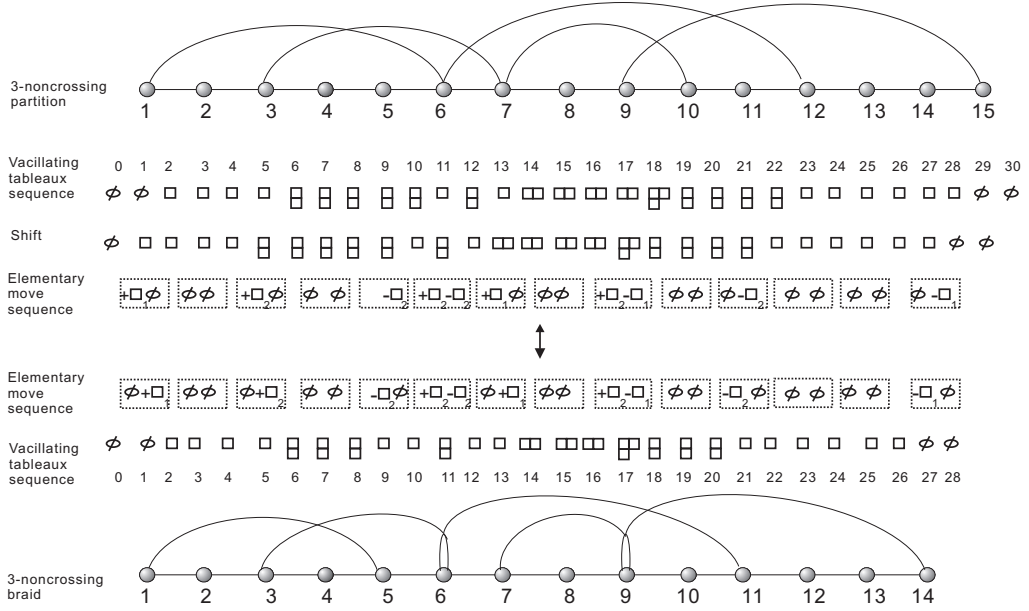,width=0.9\textwidth}\hskip15pt
 }
\caption{\small Mapping $2$-regular $3$-noncrossing partitions into
loop-free braids: illustration of Corollary~\ref{C:bi}(b). }
\label{F:map1}
\end{figure}
\begin{proof}
Assertion (a) follows from Lemma~\ref{L:bijection-2} since a
partition is completely determined by its induced
$\mathcal{P}_{W_{k-1}}$-walk of shape $\lambda^{2n-1}=\varnothing$.
(b) follows immediately from the fact that, according to the
definition of $\varphi$ given in Lemma~\ref{L:bijection-2}, any pair
of consecutive EMs $(\varnothing,
\square_1),(-\square_1,\varnothing)$ induces an EM
$(+\square_1,-\square_1)$. Therefore, $\vartheta$ maps $2$-regular,
$k$-noncrossing partitions into $k$-noncrossing braids without
loops.
\end{proof}
According to Corollary~\ref{C:bi}, each $2$-regular, $3$-noncrossing
partition on $[N]$ corresponds to a $3$-noncrossing braid without
loops over $[N-1]$. Each such braid uniquely corresponds to a
$\mathcal{B}_{W_2}$-walk of length $2N-2$ in which there does not
exist any odd $+\mathbf{e_{1}}$-step followed by an even
$-\mathbf{e_{1}}$-step. We call such a walk a
$\mathcal{B}_{W_2}^*$-walk.

\begin{figure}[ht]
\centerline{%
\epsfig{file=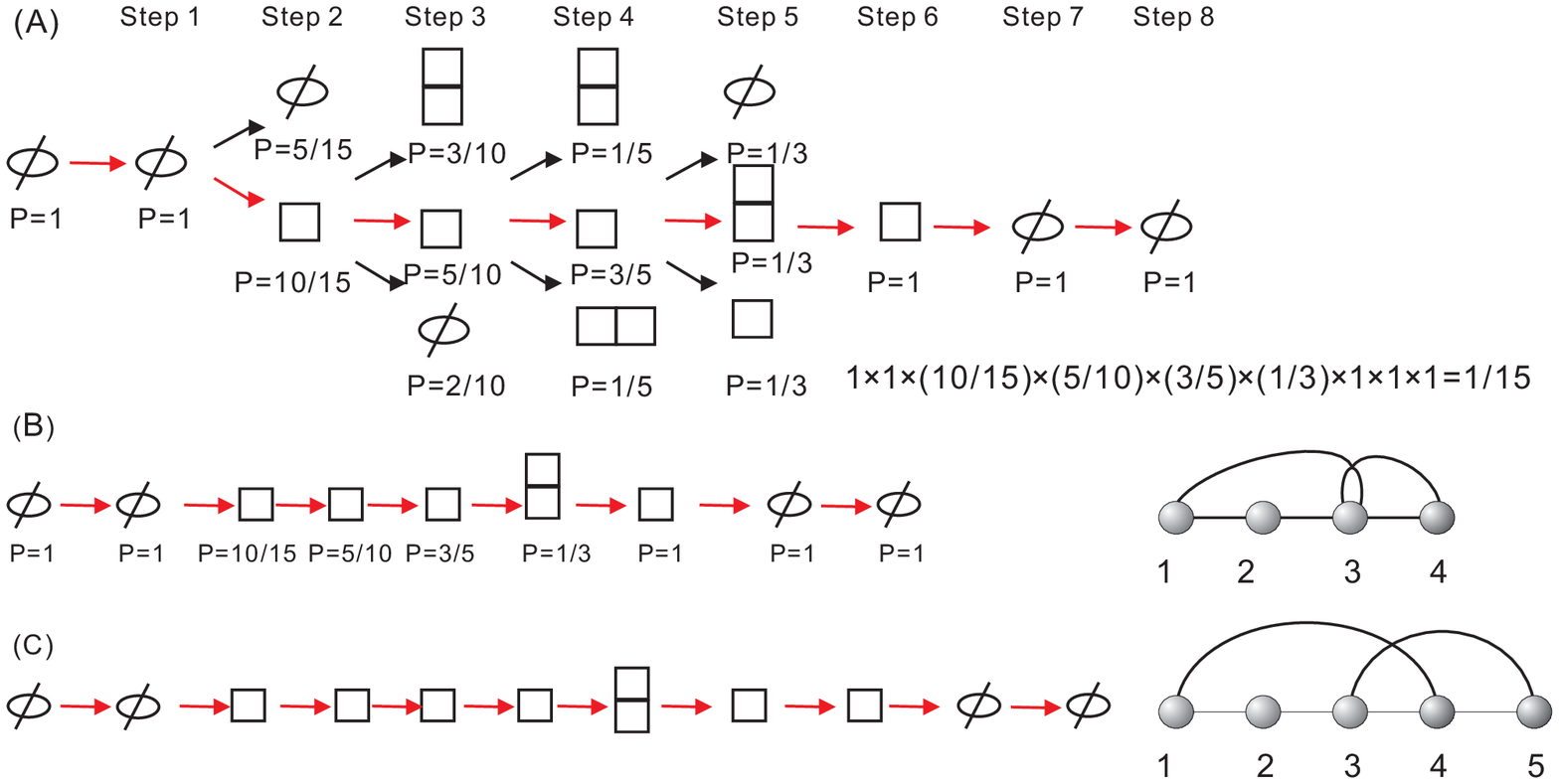,width=0.7\textwidth}\hskip15pt
 }
\caption{\small Uniform generation of $2$-regular, $3$-noncrossing
partitions: the stochastic process over shapes \textsf{(A)}, extract
a specific sampling path \textsf{ (B)} which we transform in
\textsf{(C)} into the corresponding $2$-regular, $3$-noncrossing
partition. The probabilities specified in \textsf{ (B)} are the
transition probabilities.} \label{F:1-reg}
\end{figure}

\begin{lemma}\label{L:regq}
The number of $\mathcal{B}_{W_2}^*$-walks ending at $(i,j)$ of
length $2\ell$, is given by
\begin{equation}
\sigma^{i,j;*}_{2\ell} = \sum_{h}(-1)^{h}{\ell\choose
h}\omega^{i,j}_{2(\ell-h)+1}.
\end{equation}
Furthermore, we obtain the recurrence of $\sigma^{i,j;*}_{2\ell+1}$
given by
\begin{equation}
\sigma^{i,j;*}_{2\ell+1} =
\sigma^{i-1,j;*}_{2\ell}+\sigma^{i,j-1;*}_{2\ell}+
\sigma^{i,j;*}_{2\ell}.
\end{equation}
\end{lemma}
\begin{proof}
According to Lemma~\ref{L:bijection-2}, any $\mathcal{P}_{W_2}$-walk
of length $2\ell+1$ corresponds to a unique $\mathcal{B}_{W_2}$-walk
of length $2\ell$, whence
$\sigma^{i,j}_{2\ell}=\omega^{i,j}_{2\ell+1}$. Let $A_{2\ell}(h)$
denote the set of $\mathcal{B}_{W_2}$-walks of length $2\ell$ in
which there exist at least $h$ pairs of shapes
$(\mu_{2q-1},\mu_{2q})$ induced by the EM $(+\mathbf{e_{1}},
-\mathbf{e_{1}})$. Since the removal of \NEW{$h$ pairs} of
$(+\mathbf{e_{1}}, -\mathbf{e_{1}})$-EMs from such a
$\mathcal{B}_{W_2}$-walk, results in a $\mathcal{B}_{W_2}$-walk of
length $2(\ell-h)$, we derive $A_{2\ell}(h)= {\ell \choose
h}\,\sigma^{i,j}_{2\ell-2h}$. Using the inclusion-exclusion
principle, we arrive at
\begin{equation}\label{E:IE}
\sigma^{i,j;*}_{2\ell}=\sum_{h}(-1)^{h}{\ell \choose
h}\sigma^{i,j}_{2(\ell-h)}= \sum_{h}(-1)^{h}{\ell \choose
h}\omega^{i,j}_{2(\ell-h)+1}.
\end{equation}
By construction, an odd step in a $\mathcal{B}_{W_2}$-walk is either
$-\mathbf{e_{1}}$, -$\mathbf{e_{2}}$ or $\mathbf{0}$, whence
\begin{equation*}
\sigma^{i,j;*}_{2\ell+1}=
\sigma^{i-1,j;*}_{2\ell}+\sigma^{i,j-1;*}_{2\ell}+
\sigma^{i,j;*}_{2\ell}.
\end{equation*}
\end{proof}
Via Lemma~\ref{L:regq}, we have explicit knowledge about the numbers
of $\mathcal{B}_{W_2}^*$-walks, i.e.~$\sigma^{i,j;*}_{\ell}$ for all
$i,j\in \mathbb{Z}^2$ and $1\leq \ell\leq 2N$. Accordingly, we are
now in position to generate $2$-regular, $3$-noncrossing partitions
with uniform probability via $3$-noncrossing, loop-free braids.
\begin{theorem}\label{T:main}
A random $2$-regular, $3$-noncrossing partition can be generated, in
$O(N^4)$ pre-processing time and $O(N^{3})$ space complexity, with
uniform probability in linear time. Each $2$-regular,
$3$-noncrossing partition is generated with $O(N)$ space and time
complexity.
\end{theorem}
\begin{algorithm}\label{A:algor2}
\begin{algorithmic}[1]
\STATE{\it Tableaux} (Initialize the sequence of shapes to be a list
$\{\lambda^{i}\}_{i=0}^{i=2N}$)

\STATE {$\lambda^{0}=\varnothing$, $\lambda^{2N}=\varnothing$, i=1}

\WHILE {$i < 2N$ }

\STATE {flag0 $\leftarrow$ 1, flag1 $\leftarrow$ 1, flag2
$\leftarrow$ 1, flag3 $\leftarrow$ 1}

\IF {$i$ is even}

\STATE  X[0]$\leftarrow$ ${\sf V^{*}}(\lambda^{i+1}_{\varnothing},
2n-(i+1))$

\STATE X[1]$\leftarrow$ ${\sf V^*}(\lambda^{i+1}_{+\square_{1}},
2N-(i+1))$, X[2]$\leftarrow$ ${\sf
V^{*}}(\lambda^{i+1}_{+\square_{2}}, 2n-(i+1))$

\STATE  X[3]$\leftarrow$ ${\sf V^{*}}(\lambda^{i+1}_{-\square_{1}},
2N-(i+1))$, X[4]$\leftarrow$ ${\sf
V^{*}}(\lambda^{i+1}_{-\square_{2}}, 2N-(i+1))$

\ENDIF

\IF {$i$ is odd}

\STATE  X[0]$\leftarrow$ ${\sf W^{*}}(\lambda^{i+1}_{\varnothing},
2N-(i+1))$

\STATE X[1]$\leftarrow$ ${\sf W^*}(\lambda^{i+1}_{+\square_{1}},
2N-(i+1))$, X[2]$\leftarrow$ ${\sf
W^{*}}(\lambda^{i+1}_{+\square_{2}}, 2N-(i+1))$

\STATE  X[3]$\leftarrow$ ${\sf W^{*}}(\lambda^{i+1}_{-\square_{1}},
2N-(i+1))$, X[4]$\leftarrow$ ${\sf
W^{*}}(\lambda^{i+1}_{-\square_{2}}, 2N-(i+1))$

\IF {flag0 =0}

\STATE {X[3]$\leftarrow$ 0, X[4]$\leftarrow$ 0}

\ENDIF

\IF {flag1 =0}

\STATE {X[1]$\leftarrow$ 0, X[2]$\leftarrow$ 0, X[0]$\leftarrow$ 0,
X[3]$\leftarrow$ 0}

\ENDIF

\IF {flag2 =0}

\STATE {X[1]$\leftarrow$ 0, X[2]$\leftarrow$ 0, X[3]$\leftarrow$ 0,
X[4]$\leftarrow$ 0}

\ENDIF

\IF {flag3 =0}

\STATE {X[1]$\leftarrow$ 0, X[2]$\leftarrow$ 0, X[0]$\leftarrow$ 0}

\ENDIF

\STATE {\it sum} $\leftarrow$ X[0]+X[1]+X[2]+X[3]+X[4]

\ENDIF

\STATE {\it Shape} $\leftarrow$ {\sf Random}({\it sum}) ({\sf
Random} generates the random shape $\lambda^{i+1}_{+\square_{j}(or
-\square_{j})}$ or  $\lambda^{i+1}_{\varnothing}$.)

\IF {$i$ is even and Shape=$\lambda^{i+1}_{+\square_{1}}$}

\STATE flag1 $\leftarrow$ 0

\ENDIF

\IF {$i$ is even and Shape=$\lambda^{i+1}_{-\square_{1}}$}

\STATE flag2 $\leftarrow$ 0

\ENDIF

\IF {$i$ is even and Shape=$\lambda^{i+1}_{-\square_{2}}$}

\STATE flag2 $\leftarrow$ 0

\ENDIF

\IF {$i$ is even and Shape=$\lambda^{i+1}_{+\square_{2}}$}

\STATE flag3 $\leftarrow$ 0

\ENDIF

\IF {$i$ is even and Shape=$\lambda^{i+1}_{\varnothing}$}

\STATE flag0 $\leftarrow$ 0

\ENDIF

\STATE Insert {\it Shape} into {\it Tableaux} (the sequence of
shapes)

\STATE $i \leftarrow i+1$

\ENDWHILE

\STATE  {\sf Map}({\it Tableaux}) (maps {\it Tableaux} into its
corresponding $2$-regular, $3$-noncrossing partition)
\end{algorithmic}
\caption{\small Uniform generation of $2$-regular, $3$-noncrossing
partitions.}
\end{algorithm}
\begin{proof}
We interpret $\mathcal{B}_{W_2}^*$-walks as sampling paths of a
stochastic process. To this end, we again label the $(i+1)$-th
shape, $\lambda^{i+1}_{\alpha}$ by $\alpha=\lambda^{i+1}\setminus
\lambda^{i}\in \{+\square_{1}, +\square_{2}, -\square_{1},
-\square_{2}, \varnothing\}$, where the labeling specifies the
transition from $\lambda^i$ to $\lambda^{i+1}$. In the following we
distinguish even and odd
labeled shape. \\
For $i=2s$, we let $W^{*}(\lambda^{2s}_\alpha, 2N-2s)$ denote the
number of $\mathcal{B}_{W_2}^*$-walks such that
$\alpha=\lambda^{2s}\setminus \lambda^{2s-1}$. By construction,
assume $\lambda^{2s}$ has \NEW{$(a-1)$} and $b$ boxes in its first
and second row, respectively, then we have
\begin{equation}\label{E:ws}
W^{*}(\lambda^{2s}_\alpha, 2N-2s)=\sigma^{a,b;*}_{2N-2s}
\end{equation}
i.e.~$W^{*}(\lambda^{2s}_\alpha, 2N-2s)$  is independent of $\alpha$
and we write $W^{*}(\lambda^{2s}_\alpha, 2N-2s)=
W^{*}(\lambda^{2s}, 2N-2s)$. \\
For $i=2s+1$, let $V^{*}(\lambda^{2s+1}_{\alpha},2N-(2s+1))$ denote
the number of $\mathcal{B}_{W_2}^*$-walks of shape $\lambda^{2s+1}$
of length $2N-(2s+1)$ where
$\lambda^{2s+1}\setminus\lambda^{2s}=\alpha$. Then we have setting
$u=2N-2s-2$
\begin{equation}
V^{*}(\lambda^{2s+1}_{\alpha},u+1)=
\begin{cases}
W^{*} (\lambda^{2s+2}_{+\square_{2}},u) & \text{\rm for} \
\alpha=+\square_{1}\\
W^{*}(\lambda^{2s+2}_{-\square_{1}},u)+W^{*}
(\lambda^{2s+2}_{-\square_{2}},u)
& \text{\rm for} \ \alpha=+\square_{2}\\
W^{*}(\lambda^{2s+2}_{+\square_{1}},u)+
W^{*}(\lambda^{2s+2}_{+\square_{2}},u) +
W^{*}(\lambda^{2s+2}_{\varnothing},u)
& \text{\rm for} \ \alpha=\varnothing \\
W^{*}(\lambda^{2s+2}_{\varnothing},u)
& \text{\rm for} \ \alpha= -\square_{1}, -\square_{2}.
\end{cases}
\end{equation}
We are now in position to specify the process $(X^{i})_{i=0}^{i=2N}$:\\
$\bullet$ $X^{0}=X^{2N}=\varnothing$ and $X^{i}$ is a shape with at
most 2 rows.\\
$\bullet$ for $1\leq i\leq N-1$, $(X^{2i+1}\setminus X^{2i},
X^{2i+2}\setminus
X^{2i})\in\{(-\square,\varnothing),(\varnothing,+\square),
(\varnothing,\varnothing),(+\square,-\square)\}$.\\
$\bullet$ there does not exist any subsequence $(X^{2i}, X^{2i+1},
X^{2i+2})$ such that $(X^{2i+1}\setminus X^{2i},
X^{2i+2}\setminus X^{2i})=(+\square_1,-\square_{1})$.\\
$\bullet$
the transition probabilities are given as follows:\\
{\sf(1)} for $i=2\ell$, we obtain
\begin{equation*}
\mathbb{P}(X^{i+1}=\lambda^{i+1}_{\alpha}\vert X^{i}=\lambda^{i})=
\frac{V^{*}(\lambda^{i+1}_{\alpha},2n-i-1)}{W^{*}(\lambda^{i+1},2N-i)}.
\end{equation*}
{\sf(2)} for $i=2\ell+1$, we have
\begin{equation*}
\mathbb{P}(X^{i+1}=\lambda^{i+1}_{}\vert
X^{i}=\lambda^{i}_{\alpha})=
\frac{W^{*}(\lambda^{i+1}_{},2n-i-1)}{V^{*}(\lambda^{i+1}_{\alpha},2N-i)}.
\end{equation*}
By construction,
\begin{equation}
\prod_{i=0}^{2N-1}\mathbb{P}(X^{i+1}=\lambda^{i+1}\mid
X^i=\lambda^i)= \frac{{\sf W^*}(\lambda^{2N}=\varnothing,0)} {{\sf
W^*}(\lambda^{0}=\varnothing,n)} =\frac{1}{{\sf
W^*}(\varnothing,2N)},
\end{equation}
whence $(X^i)_{i=0}^{2N}$ generates random $2$-regular,
$3$-noncrossing partitions with uniform probability in $O(N)$ time
and space, see Figure \ref{F:1-reg}.\\
\NEW{As for the derivation of the transition probabilities, supposed
that the terms $\omega^{i,j}_{\ell}$ for $0\leq i,j\leq \ell$ and
$1\leq \ell\leq 2N$ can be calculated in $O(N^4)$ time and $O(N^3)$
space
according to Theorem \ref{T:algpar}.\\
We claim for fix the indices $i_1$, $j_1$ and $s$ such that $0\leq
i_1,j_1\leq s$ and $1\leq s\leq 2N$, $\sigma^{i_1,j_1;*}_{s}$ can be
computed in $O(N)$ time. Consider the parity of $s$, there are two
cases. First, in case of $s=2\ell_1$, via using a {\bf For}-loop
summing over the terms $(-1)^{h}{\ell_1-h \choose
h}\omega^{i_1,j_1}_{2(\ell_1-h)+1}$, we obtain
$\sigma^{i_1,j_1;*}_{2\ell_1}$. Otherwise, in case of $s=2\ell_1+1$,
we first calculate $\sigma^{i-1,j_1;*}_{2\ell_1}$,
$\sigma^{i_1,j_1-1;*}_{2\ell_1}$, and $\sigma^{i_1,j_1;*}_{2\ell_1}$
via using a {\bf For}-loop according to eq.~(\ref{E:IE}) then
$\sigma^{i_1,j_1;*}_{2\ell_1+1}$ follows from
$\sigma^{i_1,j_1;*}_{2\ell_1+1}=\sigma^{i-1,j_1;*}_{2\ell_1}+
\sigma^{i_1,j_1-1;*}_{2\ell_1}+ \sigma^{i_1,j_1;*}_{2\ell_1}$.\\
Furthermore, via using three nested {\bf For} loops for $s$, $i_1$
and $j_1$ from outside to inside, we derive $\sigma^{i_1,j_1;*}_{s}$
for arbitrary $i_1$, $j_1$ and $s$. Consequently, we compute
$\sigma^{i_1,j_1;*}_{s}$ for all $i_1$, $j_1$, $s$ with
$O(N^4)+O(N^3\times N)=O(N^4)$ time and $O(N^{3})$ space complexity.
Obviously, once the terms $\sigma^{i_1,j_1;*}_{s}$ for $i_1$, $j_1$
and $s$ are calculated, we can compute the transition probabilities
in $O(1)$ time. Therefore we obtain the transition probabilities
${W^{*}(\lambda^{i+1},2N-i)}$ in $O(N^4)$ time and $O(N^{3})$ space
complexity.}
\end{proof}
\bibliographystyle{plain}
\bibliography{uni1}
\end{document}